# Improved Exponential Estimator for Population Variance Using Two Auxiliary Variables

Rajesh Singh
Department of Statistics, Banaras Hindu University(U.P.), India
(rsinghstat@yahoo.com)

Pankaj Chauhan and Nirmala Sawan
School of Statistics, DAVV, Indore (M.P.), India

Florentin Smarandache
Department of Mathematics, University of New Mexico, Gallup, USA
(smarand@unm.edu)

**Abstract**

In this paper exponential ratio and exponential product type estimators using two auxiliary variables are proposed for estimating unknown population variance $S_y^2$. Problem is extended to the case of two-phase sampling. Theoretical results are supported by an empirical study.

Key words: Auxiliary information, exponential estimator, mean squared error.

## 1. Introduction

It is common practice to use the auxiliary variable for improving the precision of the estimate of a parameter. Out of many ratio and product methods of estimation are good examples in this context. When the correlation between the study variate and the auxiliary variate is positive (high) ratio method of estimation is quite effective. On the other hand, when this correlation is negative (high) product method of estimation can be employed effectively. Let y and (x,z) denotes the study variate and auxiliary variates taking the values $y_i$ and $(x_i, z_i)$ respectively, on the unit $U_i$ (i=1,2,……,N), where x is positively correlated with y and z is negatively correlated with y. To estimate $S_y^2 = \dfrac{1}{(N-1)} \sum_{i=1}^{N} (y_i - \bar{y})^2$, it is assumed that $S_x^2 = \frac{1}{(N-1)} \sum_{i=1}^{N} (x_i - \overline{X})^2$ and $S_z^2 = \frac{1}{(N-1)} \sum_{i=1}^{N} (z_i - \overline{Z})^2$ are known. Assume that population size N is large so that the finite population correction terms are ignored.

Assume that a simple random sample of size n is drawn without replacement (SRSWOR) from U. The usual unbiased estimator of $S_y^2$ is

$$s_y^2 = \frac{1}{(n-1)} \sum_{i=1}^{n} (y_i - \overline{y})^2 \qquad (1.1)$$

where $\overline{y} = \frac{1}{n} \sum_{i=1}^{n} y_i$ is the sample mean of y.

When the population variance $S_x^2 = \frac{1}{(N-1)} \sum_{i=1}^{N} (x_i - \overline{X})^2$ is known, Isaki (1983) proposed a ratio estimator for $S_y^2$ as

$$t_k = s_y^2 \frac{S_x^2}{s_x^2} \qquad (1.2)$$

where $s_x^2 = \frac{1}{(n-1)} \sum_{i=1}^{n} (x_i - \overline{X})^2$ is an unbiased estimator of $S_x^2$.

Upto the first order of approximation, the variance of $S_y^2$ and MSE of $t_k$ (ignoring the finite population correction (fpc) term) are respectively given by

$$\text{var}(s_y^2) = \left(\frac{S_y^4}{n}\right)[\partial_{400} - 1] \qquad (1.3)$$

$$\text{MSE}(t_k) = \left(\frac{S_y^4}{n}\right)[\partial_{400} + \partial_{040} - 2\partial_{220}] \qquad (1.4)$$

where $\delta_{pqr} = \frac{\mu_{pqr}}{\left(\mu_{200}^{p/2} \mu_{020}^{q/2} \mu_{002}^{r/2}\right)}$,

$\mu_{pqr} = \frac{1}{N} \sum_{i=1}^{N} (y_i - \overline{Y})^p (x_i - \overline{X})^q (z_i - \overline{Z})^r$ ; p, q, r being the non-negative integers.

Following Bahl and Tuteja (1991), we propose exponential ratio type and exponential product type estimators for estimating population variance $S_y^2$ as –

$$t_1 = s_y^2 \exp\left[\frac{S_x^2 - s_x^2}{S_x^2 + s_x^2}\right] \tag{1.5}$$

$$t_2 = s_y^2 \exp\left[\frac{s_z^2 - S_z^2}{s_z^2 + S_z^2}\right] \tag{1.6}$$

2. **Bias and MSE of proposed estimators**

To obtain the bias and MSE of $t_1$, we write

$$s_y^2 = S_y^2(1+e_0), \quad s_x^2 = S_x^2(1+e_1)$$

Such that $E(e_0) = E(e_1) = 0$

and $E(e_0^2) = \frac{1}{n}(\partial_{400} - 1)$, $E(e_1^2) = \frac{1}{n}(\partial_{040} - 1)$, $E(e_0 e_1) = \frac{1}{n}(\partial_{220} - 1)$.

After simplification we get the bias and MSE of $t_1$ as

$$B(t_1) \cong \frac{S_y^2}{n}\left[\frac{\partial_{040}}{8} - \frac{\partial_{220}}{2} + \frac{3}{8}\right] \tag{2.1}$$

$$MSE(t_1) \cong \frac{S_y^2}{n}\left[\partial_{400} + \frac{\partial_{040}}{4} - \partial_{220} + \frac{1}{4}\right] \tag{2.2}$$

To obtain the bias and MSE of $t_2$, we write

$$s_y^2 = S_y^2(1+e_0), \quad s_z^2 = S_z^2(1+e_2)$$

Such that $E(e_0) = E(e_2) = 0$

$$E(e_2^2) = \frac{1}{n}(\partial_{004} - 1), \quad E(e_0 e_2) = \frac{1}{n}(\partial_{202} - 1)$$

After simplification we get the bias and MSE of $t_2$ as

$$B(t_2) \cong \frac{S_y^2}{n}\left[\frac{\partial_{004}}{8} + \frac{\partial_{202}}{2} - \frac{5}{8}\right] \tag{2.3}$$

$$MSE(t_2) \cong \frac{S_y^2}{n}\left[\partial_{400} + \frac{\partial_{004}}{4} + \partial_{202} - \frac{9}{4}\right] \tag{2.4}$$

3. **Improved Estimator**

Following Kadilar and Cingi (2006) and Singh et. al. (2007), we propose an improved estimator for estimating population variance $S_y^2$ as-

$$t = s_y^2 \left[ \alpha \exp\left\{ \frac{S_x^2 - s_x^2}{S_x^2 + s_x^2} \right\} + (1-\alpha) \exp\left\{ \frac{s_z^2 - S_z^2}{s_x^2 + S_z^2} \right\} \right] \tag{3.1}$$

where $\alpha$ is a real constant to be determined such that the MSE of t is minimum.

Expressing t in terms of e's, we have

$$t = S_y^2 (1+e_0) \left[ \alpha \exp\left\{ -\frac{e_1}{2}\left(1 + \frac{e_1}{2}\right)^{-1} \right\} + (1-\alpha) \exp\left\{ \frac{e_2}{2}\left(1 + \frac{e_2}{2}\right)^{-1} \right\} \right] \tag{3.2}$$

Expanding the right hand side of (3.2) and retaining terms up to second power of e's, we have

$$t \cong S_y^2 \left[ 1 + e_0 + \frac{e_2}{2} + \frac{e_2^2}{8} + \frac{e_0 e_2}{2} + \alpha\left(-\frac{e_1}{2} + \frac{e_1^2}{8}\right) - \alpha\left(\frac{e_2}{2} + \frac{e_2^2}{8}\right) \right.$$

$$\left. + e_0 \alpha\left(-\frac{e_1}{2} + \frac{e_1^2}{8}\right) - \alpha e_0 \left(\frac{e_2}{2} + \frac{e_2^2}{8}\right) \right] \tag{3.3}$$

Taking expectations of both sides of (3.3) and then subtracting $S_y^2$ from both sides, we get the bias of the estimator t, up to the first order of approximation, as

$$B(t) = \frac{S_y^2}{n} \left[ \frac{\alpha}{8}(\partial_{040} - 1) + \frac{(1-\alpha)}{8}(\partial_{004} - 1) + \frac{(1-\alpha)}{2}(\partial_{202} - 1) \right.$$

$$\left. - \frac{\alpha}{2}(\partial_{220} - 1) \right] \tag{3.4}$$

From (3.4), we have

$$(t - S_y^2) \cong S_y^2 \left[ e_0 - \frac{\alpha e_1}{2} + \frac{(1-\alpha)}{2} e_2 \right] \tag{3.5}$$

Squaring both the sides of (3.5) and then taking expectation, we get MSE of the estimator t, up to the first order of approximation, as

$$MSE(t) \cong \frac{S_y^4}{n} \left[ (\partial_{400} - 1) + \frac{\alpha^2}{4}(\partial_{040} - 1) + \frac{(1-\alpha)^2}{4}(\partial_{004} - 1) \right.$$

$$\left. - \alpha(\partial_{220} - 1) + (1-\alpha)(\partial_{202} - 1) - \frac{\alpha(1-\alpha)}{2}(\partial_{022} - 1) \right] \tag{3.6}$$

Minimization of (3.6) with respect to $\alpha$ yields its optimum value as

$$\alpha = \frac{\{\partial_{004} + 2(\partial_{220} + \partial_{202}) + \partial_{022} - 6\}}{(\partial_{040} + \partial_{004} + 2\partial_{022} - 4)} = \alpha_0 \text{(say)} \tag{3.7}$$

Substitution of $\alpha_0$ from (3.7) into (3.6) gives minimum value of MSE of t.

## 4. Proposed estimators in two-phase sampling

In certain practical situations when $S_x^2$ is not known a priori, the technique of two-phase or double sampling is used. This scheme requires collection of information on x and z the first phase sample s' of size n' (n'<N) and on y for the second phase sample s of size n (n<n') from the first phase sample.

The estimators $t_1$, $t_2$ and t in two-phase sampling will take the following form, respectively

$$t_{1d} = s_y^2 \exp\left[\frac{s_x'^2 - s_x^2}{s_x'^2 + s_x^2}\right] \qquad (4.1)$$

$$t_{2d} = s_y^2 \exp\left[\frac{s_z'^2 - s_z^2}{s_z'^2 + s_z^2}\right] \qquad (4.2)$$

$$t_d = s_y^2 \left[k \exp\left\{\frac{s_x'^2 - s_x^2}{s_x'^2 + s_x^2}\right\} + (1-k) \exp\left\{\frac{s_z'^2 - s_z^2}{s_z'^2 + s_z^2}\right\}\right] \qquad (4.3)$$

To obtain the bias and MSE of $t_{1d}$, $t_{2d}$, $t_d$, we write

$$s_y^2 = S_y^2(1+e_0), \quad s_x^2 = S_x^2(1+e_1), \quad s_x'^2 = S_x^2(1+e_1')$$

$$s_z^2 = S_z^2(1+e_2), \quad s_z'^2 = S_z^2(1+e_2')$$

where $s_x'^2 = \frac{1}{(n'-1)}\sum_{i=1}^{n'}(x_i - \bar{x}')^2$, $s_z^2 = \frac{1}{(n'-1)}\sum_{i=1}^{n'}(z_i - \bar{z}')^2$

$$\bar{x}' = \frac{1}{n'}\sum_{i=1}^{n'} x_i, \quad \bar{z}' = \frac{1}{n'}\sum_{i=1}^{n'} z_i$$

Also,

$E(e_1') = E(e_2') = 0$,

$E(e_1'^2) = \frac{1}{n'}(\partial_{040} - 1)$, $E(e_2'^2) = \frac{1}{n}(\partial_{004} - 1)$,

$E(e_1' e_2') = \frac{1}{n'}(\partial_{220} - 1)$

Expressing $t_{1d}$, $t_{2d}$, and $t_d$ in terms of e's and following the procedure explained in section 2 and section3 we get the MSE of these estimators, respectively as-

$$\text{MSE}(t_{1d}) \cong S_y^4\left[\frac{1}{n}(\partial_{400}-1)+\frac{1}{4}\left(\frac{1}{n}-\frac{1}{n'}\right)(\partial_{040}-1)\right.$$
$$\left.+\left(\frac{1}{n'}-\frac{1}{n}\right)(\partial_{220}-1)\right] \quad (4.4)$$

$$\text{MSE}(t_{2d}) \cong S_y^4\left[\frac{1}{n}(\partial_{400}-1)+\frac{1}{4}\left(\frac{1}{n}-\frac{1}{n'}\right)(\partial_{004}-1)\right.$$
$$\left.-\left(\frac{1}{n'}-\frac{1}{n}\right)(\partial_{202}-1)\right] \quad (4.5)$$

$$\text{MSE}(t_d) \cong S_y^4\left[\frac{1}{n}(\partial_{400}-1)+\frac{k^2}{4}\left(\frac{1}{n}-\frac{1}{n'}\right)(\partial_{040}-1)+\frac{(k^2-1)}{4}\left(\frac{1}{n}-\frac{1}{n'}\right)(\partial_{004}-1)\right.$$
$$+k\left(\frac{1}{n}-\frac{1}{n'}\right)(\partial_{220}-1)+(k-1)\left(\frac{1}{n'}-\frac{1}{n}\right)(\partial_{202}-1)$$
$$\left.-\frac{k(k-1)}{2}\left(\frac{1}{n'}-\frac{1}{n}\right)(\partial_{022}-1)\right] \quad (4.6)$$

Minimization of (4.6) with respect to k yields its optimum value as

$$k = \frac{\{\partial_{004}+2(\partial_{220}-1)+\partial_{022}-6\}}{(\partial_{040}+\partial_{004}+2\partial_{022}-4)} = k_0 \text{ (say)} \quad (4.7)$$

Substitution of $k_0$ from (4.7) to (4.6) gives minimum value of MSE of $t_d$.

## 5. Empirical Study

To illustrate the performance of various estimators of $S_y^2$, we consider the data given in Murthy(1967, p.-226). The variates are:

y: output, x: number of workers, z: fixed capital,

N=80, n'=25, n=10.

$\partial_{400} = 2.2667$, $\partial_{040} = 3.65$, $\partial_{004} = 2.8664$, $\partial_{220} = 2.3377$, $\partial_{202} = 2.2208$, $\partial_{400} = 3.14$

The percent relative efficiency (PRE) of various estimators of $S_y^2$ with respect to conventional estimator $s_y^2$ has been computed and displayed in table 5.1.

**Table 5.1 :  PRE of $s_y^2$, t₁, t₂ and min. MSE (t) with respect to $s_y^2$**

| Estimator | PRE$(.,s_y^2)$ |
|---|---|
| $s_y^2$ | 100 |
| t₁ | 214.35 |
| t₂ | 42.90 |
| t | 215.47 |

In table 5.2 PRE of various estimators of $s_y^2$ in two-phase sampling with respect to $S_y^2$ are displayed.

**Table 5.2 : PRE of $s_y^2$, t₁d, t₂d and min.MSE (t_d) with respect to $s_y^2$**

| Estimator | PRE $(.,s_y^2)$ |
|---|---|
| $s_y^2$ | 100 |
| t₁d | 1470.76 |
| t₂d | 513.86 |
| t_d | 1472.77 |

### 6.  Conclusion

From table 5.1 and 5.2, we infer that the proposed estimators t performs better than conventional estimator $s_y^2$ and other mentioned estimators.